\documentclass[preprint,12pt]{elsarticle}

\newtheorem{theorem}{Theorem}[section]
\newtheorem{corollary}{Corollary}[section]
\newtheorem{lemma}{Lemma}[section]
\newtheorem{remark}{Remark}[section]

\newcommand\X{{\sf X}}
\newcommand\Y{{\sf Y}}
\newcommand{\1}{{\mathbf{1}}}
\usepackage{color}
\usepackage{amssymb}
\begin{document}

\begin{frontmatter}

\title{Multistable L\'{e}vy motions and their continuous approximations}
\author{Xiequan Fan\ \ \ \ \ \ Jacques L\'{e}vy V\'{e}hel$\,^*$}
 \cortext[cor1]{\noindent Corresponding author. \\
\mbox{\ \ \ \ }\textit{E-mail}: fanxiequan@hotmail.com (X. Fan), \ \ \ \ \  jacques.levy-vehel@inria.fr (J. L\'{e}vy V\'{e}hel). }
\address{Regularity Team, Inria and MAS Laboratory, Ecole Centrale Paris - Grande Voie des Vignes,\\ 92295 Ch\^{a}tenay-Malabry, France}

\begin{abstract}
Multistable L\'{e}vy motions are extensions of L\'{e}vy motions where the stability index is allowed to vary in time. Several constructions
of these processes have been introduced recently, based on Poisson and Ferguson-Klass-LePage series representations and on multistable
measures. In this work, we prove a functional central limit theorem for  the independent-increments
multistable L\'{e}vy motion, as well as of integrals with respect to these processes, using weighted sums of independent random
variables. This allows us to construct continuous
approximations of multistable L\'{e}vy motions. In particular, we prove that multistable L\'{e}vy motions are  stochastic H\"{o}lder continuous and strongly localisable.
\end{abstract}

\begin{keyword}
(strong) localisability; multistable process; stochastic H\"{o}lder continuous; stable process; continuous approximation.
\vspace{0.3cm}
\MSC Primary  60G18, 60G17;   Secondary 60G51, 60G52.
\end{keyword}

\end{frontmatter}

% \linenumbers

\section{Introduction}
Recall that a stochastic process $\{L(t), t \geq 0\}$ is called (standard) $\alpha-$stable L\'{e}vy motion if the following three conditions hold:\\
(C1) $L(0)=0$  almost surely;\\
(C2) $L$ has independent increments;\\
(C3) $L(t)-L(s)\sim S_\alpha( (t-s)^{1/\alpha}, \beta, 0)$ for any $0\leq s < t$ and for some
$0< \alpha \leq 2, -1\leq \beta \leq 1$. Here $S_\alpha(\sigma, \beta, 0)$ stands for a stable random variable with index of
stability $\alpha$, scale parameter $\sigma$, skewness parameter $\beta$ and shift parameter
equal to 0. Recall that $\alpha$ governs the intensity of jumps.

Such processes have stationary increments, and they are $1/\alpha-$self-similar, that is, for all $c>0,$ the processes $\{L(c\,t), t\geq 0\}$ and $\{c^{1/\alpha}L(t), t\geq 0\}$ have the same finite-dimensional distributions.
An $\alpha-$stable L\'{e}vy motion is symmetric when $\beta=0$. Stable L\'{e}vy motions, and, more generally, stable processes
have been the subject of intense activity in recent years, both on the theoretical side (see, e.g. \cite{ST94}) and in
applications \cite{N12}. However, the stationary property of their increments restricts their use in some situations,
and generalizations are needed for instance to model
real-world phenomena such as financial records, epileptic episodes in EEG or internet traffic. A significant feature in these
cases is that the ``local intensity of jumps'' varies with time $t$. A way to deal with such a variation is set up a class of processes
whose stability index $\alpha$ is a function of $t$. More precisely, one aims at defining non-stationary increments processes which
are, at each time $t$, ``tangent'' (in a certain sense explained below) to a stable process with stability index $\alpha(t)$.

Formally, one says that a stochastic process $\{X(t), t \in [0, 1]\}$ is {\it multistable} \cite{FL09} if, for almost
all $t\in [0,1)$, $X$ is \emph{localisable} at $t$ with tangent process $X_t '$ an $\alpha(t)-$stable process. Recall that $\{X(t), t \in [0, 1]\}$ is said to be $h-$localisable at $t$ (cf.\ \cite{F,F2}),
with $h>0$, if there exists a non-trivial process $X_t '$, called the tangent process of $X$ at $t$, such that
\begin{eqnarray}\label{indis}
\lim_{r\searrow 0} \frac{X(t+ru)-X(t)}{r^h} = X_t'(u),
\end{eqnarray}
where convergence is in finite dimensional distributions.

Let $D[0, 1]$ be the set of c\`{a}dl\`{a}g functions on $[0, 1],$ that is functions which are continuous
on the right and have left limits at all $t \in [0, 1],$ endowed with the Skorohod  metric $d_S$ \cite{B68}.
If $X$ and $X_t'$ have versions in $D[0, 1]$ and convergence
in (\ref{indis}) is in distribution with respect to $d_S$, one says that $X$ is $h-$\emph{strongly localisable}
at $t$ with strong local form $X_t'$.

In this work, we will be concerned with the simplest non-trivial multistable processes, namely multistable L\'evy motions (MsLM), which
are non-stationary increments extensions of stable L\'{e}vy motions. Two such extensions exist \cite{FL09,FL12}:
\begin{enumerate}
\item The {\it field-based} MsLM admit the following series representation:
\begin{equation} \label{FBP}
L_{F}(t)  = C_{ \alpha(t)}^{1/ \alpha(t)}\sum_{(\X,\Y)\in\Pi}
\1_{[0,t]}(\X)\Y^{<-1/\alpha(t)>} \quad \quad  (t \in [0,T]),
\end{equation}
where $\Pi$ is a Poisson point process on $[0,1] \times \mathbb{R}$ with mean measure the Lebesgue measure $\mathcal{L}$,
$a^{<b>} := \mbox{sign}(a)|a|^{b}$ and
\begin{equation}\label{calphax}
C_{u}= \left( \int_{0}^{\infty} x^{-u} \sin (x)dx \right)^{-1}.
\end{equation}
Their joint characteristic function reads:
\begin{equation}\label{FBCF}
\mathbb{E}  \exp \left\{i \sum_{j=1}^{m}\limits \theta_j L_{F}(t_j) \right\}  = \exp \left\{ -2 \int_{[0,T]} \int_{0}^{+ \infty} \sin^2\Bigg( \sum_{j=1}^{m} \theta_j \frac{C_{\alpha(t_j)}^{1/\alpha(t_j)}}{2y^{1/\alpha(t_j)}} \mathbf{1}_{[0,t_j]}(x) \Bigg)\hspace{0.1cm} dy \hspace{0.1cm} dx \right\}
\end{equation}
for $ d \in \mathbb{N}, (\theta_1, \ldots, \theta_d) \in \mathbb{R}^d$ and $(t_1, \ldots , t_d) \in \mathbb{R}^d.$
These processes have correlated increments, and they are localisable as soon as the function $\alpha$ is
H\"{o}lder-continuous.
\item The {\it independent-increments} MsLM admit  the following series representation:
\begin{equation}\label{PRLII}
L_{I}(t)=\sum_{(\X,\Y)\in \Pi}C_{\alpha(\X)}^{1/ \alpha(\X)}\mathbf{1}_{[0,t]}(\X)\Y^{<-1/\alpha(\X)>} \quad \quad  (t \in [0,T]).
\end{equation}
As their name indicates, they have independent increments, and their joint characteristic function reads:
\begin{eqnarray}
\mathbb{E}  \exp\left\{ i \sum_{j=1}^d \theta_j L_{I }(t_j) \right\}  =\exp \left\{- \int \Big| \sum_{j=1}^d \theta_j \mathbf{1}_{[0, \ t_j] }(s)\Big|^{\alpha(s)} ds \right\}, \label{hjdns}
\end{eqnarray}
for $ d \in \mathbb{N}, (\theta_1, \ldots, \theta_d) \in \mathbb{R}^d$ and $(t_1, \ldots , t_d) \in \mathbb{R}^d.$
These  processes are  localisable as soon as the function $\alpha$ verifies:
\begin{eqnarray}\label{coalph}
\Big(\alpha (x)-\alpha (x+t) \Big)\ln t \rightarrow 0
\end{eqnarray}
uniformly for all $x$ in finite interval as $t\searrow0$ \cite{FL12}.
\end{enumerate}
Of course, when $\alpha(t)$ is a constant $\alpha$ for all $t$, both $L_{F}$ and $L_I$ are simply the Poisson representation of
$\alpha-$stable L\'evy motion, that we denote by $L_{\alpha}$. In general, $L_F$ and $L_I$ are semi-martingales \cite{GLL13b}.
For more properties of $L_F$, such as Ferguson-Klass-LePage series representations and H\"{o}lder exponents, we refer to \cite{FLV09,GLL12,GLL13a}.

In this paper, we prove a functional central limit theorem for  independent-increments  MsLM:
we show that certain weighted sums of independent random variables converge in $(D[0,1],d_S)$ to $L_I.$
This allows us to obtain strong localisability of these processes.  Moreover, we establish
continuous approximations of MsLM and an alternative representation for the integrals of multistable L\'{e}vy measure.
 Some properties of the integrals of multistable L\'{e}vy measure are investigated. In particular,  we prove that MsLM are stochastic H\"{o}lder continuous and strongly localisable.

The paper is organized as follows. In Section \ref{secmlevy}, we present the  functional central limit theorem for  independent-increments MsLM. In Section \ref{contapp},
we establish continuous approximations of MsLM. In the last section, we give a representation of MsLM and investigate some properties, including stochastic H\"{o}lder continuous and strongly localisable, of the integrals of multistable
L\'{e}vy measure.

\section{Functional Central Limit Theorems for Multistable L\'{e}vy  Motions}\label{secmlevy}
We show in this section how to approximate the independent-increments MsLM in law by weighted sums of independent random variables.
\begin{theorem}\label{fhnk} Let $(\alpha_n(u))_n,\alpha(u), u  \in [0,1],$ be a class of c\`{a}dl\`{a}g functions
ranging in $[a,b] \subset (0,2]$ such that the sequence $(\alpha)_n$ tends to $\alpha$ in the uniform metric.
Let $\big( X(k,n) \big)_{n \in \mathbb{N}, \ k=1,...,2^n} $ be a family of independent and  symmetric
$\alpha_n(\frac{k}{2^n})-$stable random variables with unit scale parameter, i.e.,
$X(k,n) \sim S_{\alpha_n(\frac{k}{2^n}) }(1, 0 , 0 )$. Then
\begin{itemize}
\item the sequence of processes
\begin{eqnarray}\label{multilevyI}
L^{(n)}_{I}(u) =   \sum_{k=1}^{\lfloor 2^n u \rfloor} \Big(\frac{1}{2^n}\Big)^{1/\alpha_n(\frac{k}{2^n})} X(k,n) , \ \ \  \ \ \ \ \ \ \ u  \in [0,1],
\end{eqnarray}
tends in distribution to $L_I(u)$ in $(D[0,1],d_S),$ where $\lfloor x \rfloor$ is the largest integer smaller than or equal to $x$.
In particular, if $\alpha$ satisfies condition (\ref{coalph}), then $L_I(u)$ is localisable at all times.
\item the sequence of processes
\begin{eqnarray}\label{multilevyII}
L^{(n)}_{R}(u)   = \sum_{k=1}^{\Gamma_{\lfloor 2^n u \rfloor}} \Big(\frac{1}{2^n}\Big)^{1/\alpha_n(\frac{k}{2^n})} X(k,n) , \ \ \  \ \ \ \ \ \ \ u  \in [0,1],
\end{eqnarray}
tends in distribution to $L_I(u)$ in $(D[0,1],d_S),$ where $(\Gamma_i)_{i\geq 1}$ is a sequence of arrival times of a Poisson process with unit arrival rate and
is independent of $\big( X(k,n) \big)_{n \in \mathbb{N}, \ k=1,...,2^n}$.
\item the sequence of processes
\begin{eqnarray}\label{multilevyC}
L^{(n)}_{C}(u)   = \sum_{k=1}^{\lfloor 2^n u \rfloor} \Big(\frac{1}{\Gamma_{2^n}}\Big)^{1/\alpha_n(\frac{k}{2^n})} X(k,n) , \ \ \  \ \ \ \ \ \ \ u  \in [0,1],
\end{eqnarray}
tends in distribution to $L_I(u)$ in $(D[0,1],d_S).$
\end{itemize}
\end{theorem}
\noindent\emph{Proof.}
We prove the first claim  by the following three steps.

First, we prove that $L^{(n)}_{I}(u) $ converges to $L_{I}(u)$ in finite dimensional distribution.   For any $u_1,u_2 \in [0,1]$ and $u_2> u_1,$ we have, for any $\theta \in \mathbb{R},$
\begin{eqnarray}
\lim_{n \rightarrow \infty} \mathbb{E} e^{i\theta \big(L^{(n)}_{I}( u_2) - L^{(n)}_{I}(u_1) \big ) }
&=& \lim_{n \rightarrow \infty}   \exp\left\{ -\sum_{k= \lfloor2^nu_{1}\rfloor+1}^{\lfloor2^nu_{2}\rfloor} \frac{1}{\ 2^n} |\theta|^{\alpha_n(\frac{k }{2^n})} \right\}. \label{chfunction}
\end{eqnarray}
Notice  that
\begin{eqnarray}
\sum_{k= \lfloor2^nu_{1}\rfloor+1}^{\lfloor2^nu_{2}\rfloor} \frac{1}{\ 2^n} \Big| |\theta|^{\alpha_n(\frac{k }{2^n})} -|\theta|^{\alpha (\frac{k }{2^n})} \Big| &\leq&   |\theta|^\tau \log |\theta|  \sum_{k= \lfloor2^nu_{1}\rfloor+1}^{\lfloor2^nu_{2}\rfloor} \frac{1}{\ 2^n}  \Big| \alpha_n\Big(\frac{k }{2^n}\Big)   -\alpha \Big(\frac{k }{2^n}\Big)   \Big|  \nonumber \\
&\leq& \Big|\Big| \alpha_n(\cdot)   -\alpha (\cdot) \Big|\Big|_{\infty}  |\theta|^\tau \log |\theta|
 \frac{\lfloor2^nu_{2}\rfloor - \lfloor2^nu_{1}\rfloor }{\ 2^n} , \label{gbsd23}
\end{eqnarray}
where $\tau=a\mathbf{1}_{[0,\, 1 )}(|\theta|)+ b\mathbf{1}_{[1,\, \infty)} (|\theta|).$ By hypothesis, we have $$ \lim_{n\rightarrow \infty}  || \alpha_n(\cdot)   -\alpha (\cdot) ||_{\infty}=0.$$ Thus  inequality (\ref{gbsd23}) implies that
\begin{eqnarray*}
\lim_{n \rightarrow \infty} \sum_{k= \lfloor2^nu_{1}\rfloor+1}^{\lfloor2^nu_{2}\rfloor} \frac{1}{\ 2^n} |\theta|^{\alpha_n(\frac{k }{2^n})}& = &
\lim_{n \rightarrow \infty} \sum_{k= \lfloor2^nu_{1}\rfloor+1}^{\lfloor2^nu_{2}\rfloor} \frac{1}{\ 2^n} |\theta|^{\alpha(\frac{k }{2^n})} \\
& = &  \int_ { u_1}^ {u_2 }   |\theta|^{\alpha (s)}  ds .
\end{eqnarray*}
From (\ref{chfunction}), it follows that
\begin{eqnarray}
\lim_{n \rightarrow \infty} \mathbb{E} e^{i\theta \big(L^{(n)}_{I}( u_2) - L^{(n)}_{I}(u_1) \big ) }
&=&   \exp\left\{ - \int_ { u_1}^ {u_2 }   |\theta|^{\alpha (s)}  ds   \right\}.  \label{hjkls}
\end{eqnarray}
Hence $ L^{(n)}_{I}( u_2) - L^{(n)}_{I}(u_1) $ converges in distribution and the characteristic function of its limit is defined by (\ref{hjkls}).  Since $L^{(n)}_{I}(u) $ has independent increments,  the limit of $L^{(n)}_{I}( u)$ has the joint characteristic function (\ref{hjdns}), i.e., $L^{(n)}_{I}( u)$ converges to
$L_{I}( u)$ in finite dimensional distribution.

Second, we prove that $L^{(n)}_{I}( u)$ converges to
$L_{I}( u)$  in $(D[0,1],d_S)$.   By Theorem 15.6 of Billingsley \cite{B68}, it suffices to show that
\begin{eqnarray}\label{sfffdf}
\mathbb{P}\Big( \Big| L^{(n)}_{I}(u)-L^{(n)}_{I}(u_1)\Big|\geq \lambda,\ \Big|L^{(n)}_{I}(u_2)-L^{(n)}_{I}(u)\Big|\geq \lambda \Big) \leq \frac{C}{\lambda^{2\gamma } } \Big[  u_2 - u_1 \Big]^2
\end{eqnarray}
for $u_1 \leq u \leq u_2, \lambda> 0$ and $n\geq 1$, where $\gamma  = a\mathbf{1}_{[2,\, \infty)} (\lambda)+ b\mathbf{1}_{(0,\, 2 )}(\lambda)$ and
$C$ is a constant depending only on $a$ and $b.$
If $u_2-u_1 < 1/2^n$, then either $L^{(n)}_{I}(u_2)=L^{(n)}_{I}(u)$ or $ L^{(n)}_{I}(u)=L^{(n)}_{I}(u_1)$; in either of 
these cases the left side of (\ref{sfffdf}) vanished.  
Next,  we  consider the case of $u_2-u_1 \geq 1/2^n.$  
Since $L^{(n)}_{I}(u)-L^{(n)}_{I}(u_1)$ and $L^{(n)}_{I}(u_2)-L^{(n)}_{I}(u)$ are independent, it follows that
\begin{eqnarray*}
&&\mathbb{P}\Big( \Big| L^{(n)}_{I}(u)-L^{(n)}_{I}(u_1)\Big|\geq \lambda,\ \Big|L^{(n)}_{I}(u_2)-L^{(n)}_{I}(u)\Big|\geq \lambda \Big)=
 \\ && \ \ \ \ \ \ \ \ \ \ \ \ \ \ \ \ \ \ \ \ \ \ \mathbb{P}\Big( \Big| L^{(n)}_{I}(u)-L^{(n)}_{I}(u_1)\Big|\geq \lambda\Big) \ \mathbb{P}\Big(\Big|L^{(n)}_{I}(u_2)-L^{(n)}_{I}(u)\Big|\geq \lambda \Big).
\end{eqnarray*}
Then, by the Billingsley inequality (cf.  p.\ 47 of \cite{B68}), it is easy to see that
\begin{eqnarray*}
\mathbb{P}\Big( \Big| L^{(n)}_{I}(u)-L^{(n)}_{I}(u_1)\Big|\geq \lambda  \Big) &\leq& \frac{\lambda}{2} \int_{-2/\lambda}^{2/\lambda } \Bigg( 1- \mathbb{E}e^{i\theta \big(L^{(n)}_{I}(u)-L^{(n)}_{I}(u_1) \big)}\Bigg) d\theta \\
&=& \frac{\lambda}{2} \int_{-2/\lambda}^{2/\lambda } \Bigg( 1-  \exp\Bigg\{- \sum_{k=\lfloor2^nu_1\rfloor+1}^{\lfloor2^nu\rfloor} \frac{1}{2^n} |\theta |^{\alpha_n(\frac{k}{2^n})} \Bigg\}\Bigg) d\theta \\
&\leq& \frac{\lambda}{2} \int_{-2/\lambda}^{2/\lambda } \sum_{k=\lfloor2^nu_1\rfloor+1}^{\lfloor2^nu\rfloor} \frac{1}{2^n} |\theta |^{\alpha_n(\frac{k}{2^n})}  d\theta \\
&\leq&\sum_{k=\lfloor2^nu_1\rfloor+1}^{\lfloor2^nu\rfloor} \frac{1}{2^n}  \frac{\lambda}{2} \int_{-2/\lambda}^{2/\lambda } \Big|  \theta\Big|^{\gamma}      \, d\theta  \\
&\leq&\frac{C_1}{\lambda^{\gamma} } \Bigg[ \frac{\lfloor2^nu\rfloor- \lfloor2^nu_1\rfloor}{2^n}    \Bigg],
\end{eqnarray*}
where  $C_1$ is  a constant  depending only on $a$ and $b.$
 Similarly, it holds
\begin{eqnarray}
\mathbb{P}\Big( \Big|L^{(n)}_{I}(u_2)-L^{(n)}_{I}(u)\Big|\geq \lambda \Big) \leq \frac{C_2}{\lambda^{\gamma} }\Bigg[ \frac{\lfloor2^n u_2\rfloor- \lfloor2^n u\rfloor}{2^n}    \Bigg] ,
\end{eqnarray}
where  $C_2$ is  a constant  depending only on $a$ and $b$.  Using the inequality $xy \leq (x+y)^2/4$ for all $x, y \geq0,$  we deduce
\begin{eqnarray*}
&& \mathbb{P}\Big( \Big| L^{(n)}_{I}(u)-L^{(n)}_{I}(u_1)\Big|\geq \lambda,\ \Big|L^{(n)}_{I}(u_2)-L^{(n)}_{I}(u)\Big|\geq \lambda \Big) \\
&\leq&  \frac{C_1C_2}{\lambda^{ 2\gamma } } \Bigg[ \frac{\lfloor2^nu\rfloor- \lfloor2^nu_1\rfloor}{2^n}    \Bigg]\Bigg[ \frac{\lfloor2^nu_2\rfloor- \lfloor2^nu\rfloor}{2^n}    \Bigg] \\
& \leq& \frac{C_1C_2}{4}  \frac{1}{\lambda^{ 2\gamma } } \Bigg[  \frac{\lfloor2^nu_2\rfloor- \lfloor2^nu_1\rfloor}{2^n}   \Bigg]^2\\
& \leq& C_1C_2  \frac{1}{\lambda^{ 2\gamma } } \Big[ u_2- u_1    \Big]^2,
\end{eqnarray*}
where the last line follows from the fact that
\begin{eqnarray*}
 \frac{\lfloor2^nu_2\rfloor- \lfloor2^nu_1\rfloor}{2^n}  
 \ \leq \ \frac{2^nu_2- 2^nu_1+ 1}{2^n} \  \leq\  2\Big[ u_2- u_1    \Big].
\end{eqnarray*}
This completes the proof of (\ref{sfffdf}).

Third, we prove that if $\alpha$ satisfies condition (\ref{coalph}), then   $L_{I}( u)$  is localisable at all times. Falconer and  Liu (cf.\ Theorem 2.7 of \cite{FL12}) have proved that the  process $L_{I}( u)$, defined by the joint characteristic function (\ref{hjdns}), is localisable  at $u$ to L\'{e}vy  motions $L_{\alpha (u)}(\cdot)$ with the stability index $\alpha (u)$.
Here we give another proof to complete our argument.  For any $(t_1, ..., t_d) \in [0,1]^d,$ from equality (\ref{hjdns}), it is easy to see that
\begin{eqnarray}
&&  \mathbb{E} \exp\left\{ i \sum_{j=1}^d \theta_j \left( \frac{L_{I}(u+rt_j)- L_{I}(u )}{ r^{1/\alpha(u)} } \right) \right\}  \nonumber\\
&&\ \ \ \ \ \ \ \ \ \ \ \ \ \ \ \ =\ \exp \left\{ - \int \Big| \sum_{j=1}^d \theta_jr^{-1/\alpha(u)} \mathbf{1}_{[u,\  u+rt_j] }(s)\Big|^{\alpha (s)} ds \right\}. \nonumber
\end{eqnarray}
Setting $s=u+rt$, we find that
\begin{eqnarray}
&&\mathbb{E}  \exp\left\{ i \sum_{j=1}^d \theta_j \left( \frac{L_{I}(u+rt_j)- L_{I}(u )}{ r^{1/\alpha(u)} } \right) \right\}  \ \ \ \ \ \ \ \ \ \ \ \ \ \ \ \ \nonumber\\
&&\ \ \ \ \ \ \ \ \ \ \ \ \ \ \ \ =\ \exp \left\{ - \int \Big| \sum_{j=1}^d \theta_j \mathbf{1}_{[0,\  t_j] }(t)\Big|^{\alpha (u+rt)} r^{(\alpha (u)-\alpha (u+rt))/\alpha (u)} dt \right\}. \nonumber
\end{eqnarray}
By condition (\ref{coalph}), it follows that
\begin{eqnarray}\label{sdfdsdf}
\lim_{r\searrow 0} r^{(\alpha (u)-\alpha (u+rt))/\alpha (u)}=1 \ \ \ \ \ \textrm{and} \ \ \ \ \  \lim_{r\searrow 0} \alpha(u+rt)=\alpha(u).
\end{eqnarray}
Hence, using dominated convergence theorem,   we have
\begin{eqnarray}
 \lim_{r\searrow 0}\mathbb{E}  \exp\left\{ i \sum_{j=1}^d \theta_j \left( \frac{L_{I}(u+rt_j)- L_{I}(u )}{ r^{1/\alpha(u)} } \right) \right\}
&=& \exp \left\{ - \int \Big| \sum_{j=1}^d \theta_j \mathbf{1}_{[0, \  t_j] }(t)\Big|^{\alpha (u)}  dt \right\} \nonumber \\
&=& \mathbb{E}  \exp\left\{ i \sum_{j=1}^d \theta_jL_{\alpha (u)}(t_j )   \right\} , \nonumber
\end{eqnarray}
which means that $L_{I}( u)$  is localisable at $u$ to an $\alpha (u)-$stable L\'{e}vy motion  $L_{\alpha (u)} (t).$
This completes the proof of the first claim of the theorem.

Next, we prove the second claim of the theorem. For any $u_1,u_2 \in [0,1]$ and $u_2> u_1,$  it is easy to see that, for any $\theta \in \mathbb{R},$
\begin{eqnarray}
\lim_{n \rightarrow \infty} \mathbb{E} e^{i\theta (L^{(n)}_{R}( u_2) - L^{(n)}_{R}(u_1)  ) }
 &=& \lim_{n \rightarrow \infty}  \mathbb{E} \exp\left\{ -\sum_{k= \Gamma_{\lfloor2^nu_1\rfloor}+1}^{\Gamma_{\lfloor2^nu_2\rfloor}} \frac{1}{\ 2^n} |\theta|^{\alpha_n(\frac{k }{2^n})} \right\} \nonumber\\
&=&   \exp\left\{ - \int_ { u_1}^ {u_2 }   |\theta|^{\alpha (s)}  ds   \right\},
\end{eqnarray}
where the last line follows from the weak law of large numbers. Notice that $L^{(n)}_{R}(u) $ also has independent increments. The rest of the proof of the second claim is similar to the proof of the first one.
For this reason, we shall not carry it out.

In the sequel, we prove the third claim by the following two steps.

First, we prove that $L^{(n)}_{C}(u) $ converges to $L_{I}(u)$ in finite dimensional distribution. It is worth noting that   $L^{(n)}_{C}$ does not have independent increments.  This property
implies that we cannot use the previous method.   For any $(u_1, ..., u_d) \in [0,1]^d$ and any $(\theta_1, ..., \theta_d) \in \mathbb{R}^d$ such that $0=u_0\leq u_1\leq u_2 \leq ...\leq u_d,$ we have
\begin{eqnarray}
\lim_{n \rightarrow \infty} \mathbb{E} e^{i \sum_{j=1}^d \theta_j  L^{(n)}_{C}( u_j)  }
&=& \lim_{n \rightarrow \infty}    \mathbb{E} \exp\left\{ i \sum_{l=1}^d \sum_{k= \lfloor 2^n   u_{l-1}\rfloor}^{\lfloor 2^n   u_l\rfloor}  \sum_{j=l}^d \theta_j  \Big(\frac{1}{\Gamma_{2^n}}\Big)^{1/\alpha_n(\frac{k}{2^n})} X(k,n)\right\} \nonumber\\
&=& \lim_{n \rightarrow \infty}  \mathbb{E} \exp\left\{ -\sum_{l=1}^d \sum_{k= \lfloor 2^n   u_{l-1}\rfloor}^{\lfloor 2^n   u_l\rfloor} \Big| \sum_{j=l}^d \theta_j\Big|^{\alpha_n(\frac{k }{2^n})}    \frac{1}{ 2^n} \frac{2^n}{\Gamma_{2^n}} \right\} \nonumber\\
&=&   \exp\left\{ -\sum_{l=1}^d \int_ {u_{l-1}}^ {u_l }   \Big|\sum_{j=l}^d \theta_j \Big|^{\alpha (s)}  ds   \right\}\nonumber\\
&=&   \exp\left\{ - \int  \Big| \sum_{j=1}^d \theta_j \mathbf{1}_{ [0, u_j) }(s) \Big|^{\alpha (s)}  ds   \right\},\nonumber
\end{eqnarray}
which gives the  joint characteristic function of $L_I.$

Second,   we prove that $L^{(n)}_{C}( u)$ converges to
$L_{I}( u)$  in $(D[0,1],d_S)$.   Again by Theorem 15.6 of Billingsley \cite{B68}, it  suffices  to show that
\begin{eqnarray}\label{sfffds}
\mathbb{P}\Big( \Big| L^{(n)}_{C}(u)-L^{(n)}_{C}(u_1)\Big|\geq \lambda,\ \Big|L^{(n)}_{C}(u_2)-L^{(n)}_{C}(u)\Big|\geq \lambda \Big) \leq \frac{C}{\lambda^{2\gamma } } \Big[  u_2 - u_1 \Big]^2
\end{eqnarray}
for $u_1 \leq u \leq u_2, \lambda> 0$ and $n\geq 1$, where $\gamma  = a\mathbf{1}_{[2,\, \infty)} (\lambda)+ b\mathbf{1}_{(0,\, 2 )}(\lambda)$ and
$C$ is a constant depending only on $a$ and $b.$  We need only consider the case of $u_2-u_1 \geq 1/2^n.$  
Since $L^{(n)}_{C}(u)-L^{(n)}_{C}(u_1)$ and $L^{(n)}_{C}(u_2)-L^{(n)}_{C}(u)$
are conditionally independent given $\Gamma_{2^n}$,
it follows that
\begin{eqnarray}
&&\mathbb{P}\Big( \Big| L^{(n)}_{C}(u)-L^{(n)}_{C}(u_1)\Big|\geq \lambda,\ \Big|L^{(n)}_{C}(u_2)-L^{(n)}_{C}(u)\Big|\geq \lambda \ \Big|\ \Gamma_{2^n} \Big)=
 \nonumber \\ && \ \ \ \ \ \ \ \ \ \ \    \mathbb{P}\Big( \Big| L^{(n)}_{C}(u)-L^{(n)}_{C}(u_1)\Big|\geq \lambda \ \Big|\ \Gamma_{2^n}\Big) \ \mathbb{P}\Big(\Big|L^{(n)}_{C}(u_2)-L^{(n)}_{C}(u)\Big|\geq \lambda \ \Big|\ \Gamma_{2^n} \Big). \label{ineq18}
\end{eqnarray}
It is easy to see that
\begin{eqnarray*}
&&\mathbb{P}\Big( \Big| L^{(n)}_{C}(u)-L^{(n)}_{C}(u_1)\Big|\geq \lambda \ \Big|\ \Gamma_{2^n} \Big) \ \leq\ \frac{\lambda}{2} \int_{-2/\lambda}^{2/\lambda } \Bigg( 1- \mathbb{E} \Big[e^{i\theta \big(L^{(n)}_{C}(u)-L^{(n)}_{C}(u_1)\big)} \ \Big|\ \Gamma_{2^n}\Big] \Bigg) d\theta \\
&& \ \ \ \  \ \ \ \ \ \ \  \ \ \  \ \ \ \  \ \ \ \ \  \ \ \  \ \ \ = \ \frac{\lambda}{2} \int_{-2/\lambda}^{2/\lambda } \Bigg( 1- \mathbb{E} \Big[ \exp\Big\{- \sum_{k=\lfloor2^nu_1\rfloor+1}^{\lfloor2^nu\rfloor} \frac{1}{\Gamma_{2^n}} |\theta |^{\alpha_n(\frac{k}{2^n})} \Big\}\ \Big|\ \Gamma_{2^n}\Big] \Bigg) d\theta \\
&& \ \ \ \  \ \ \ \ \ \ \  \ \ \  \ \ \ \  \ \ \ \ \  \ \ \  \ \ \ \leq \ \frac{\lambda}{2} \int_{-2/\lambda}^{2/\lambda }\sum_{k=\lfloor2^nu_1\rfloor+1}^{\lfloor2^nu\rfloor} \frac{1}{ \Gamma_{2^n}} |\theta |^{\alpha_n(\frac{k}{2^n})}   d\theta \\
&& \ \ \ \  \ \ \ \ \ \ \  \ \ \  \ \ \ \  \ \ \ \ \  \ \ \  \ \ \  \leq\ \frac{C_1}{\lambda^{\gamma} } \frac{2^n}{ \Gamma_{2^n}} \Bigg[ \frac{\lfloor2^nu\rfloor- \lfloor2^nu_1\rfloor}{2^n}    \Bigg],
\end{eqnarray*}
where $C_1$ is  a constant  depending only on $a$ and $b.$
 Similarly, it holds
\begin{eqnarray}
\mathbb{P}\Big( \Big|L^{(n)}_{C}(u_2)-L^{(n)}_{C}(u)\Big|\geq \lambda \ \Big|\ \Gamma_{2^n}\Big) \leq \frac{C_2}{\lambda^{\gamma} } \frac{2^n}{ \Gamma_{2^n}}\Bigg[ \frac{\lfloor2^nu_2\rfloor- \lfloor2^nu\rfloor}{2^n}    \Bigg] ,
\end{eqnarray}
where  $C_2$ is  a constant  depending only on $a$ and $b$.  From (\ref{ineq18}),   we find
\begin{eqnarray*}
&&\mathbb{P}\Big( \Big| L^{(n)}_{C}(u)-L^{(n)}_{C}(u_1)\Big|\geq \lambda,\ \Big|L^{(n)}_{C}(u_2)-L^{(n)}_{C}(u)\Big|\geq \lambda \Big) \\
&\leq& \frac{C_1C_2}{\lambda^{ 2\gamma } } \Bigg[ \frac{\lfloor2^nu\rfloor- \lfloor2^nu_1\rfloor}{2^n}    \Bigg]\Bigg[ \frac{\lfloor2^nu_2\rfloor- \lfloor2^nu\rfloor}{2^n}    \Bigg] \, \mathbb{E}\Big[ \Big(\frac{2^n}{ \Gamma_{2^n}} \Big)^2 \Big]\\
&\leq& \frac{C  }{\lambda^{ 2\gamma } }  \Big[  u_2   - u_1 \Big]^2 ,
\end{eqnarray*}
where  $C$ is  a constant  depending only on $a$ and $b$.
This completes the proof of (\ref{sfffds}). \hfill\qed

\begin{remark}
Let us comment on Theorem \ref{fhnk}.
\begin{enumerate}
\item We can define the  independent-increments  MsLM $\{L_{I}(x): \ x \in \mathbb{R} \}$ on the whole line as follows.  Let $\alpha(x), \ x \in \mathbb{R},$  be a continuous function ranging in $[a,b] \subset (0,2]$, and satisfies condition  (\ref{coalph})  uniformly for all $x$ in finite interval as $t \searrow 0$.
Set the functions $\alpha_k(x)=\alpha(x+k)$ for all $k \geq 0$ and $x \in [0, 1].$  For any $\alpha_k(x)$, by Theorem \ref{fhnk}, we can construct MsLM $$L_{I_k }(x): [0,1] \rightarrow \mathbb{R},\ \ \ \ \ k\geq 0.$$  Taking a sequence of independent processes $L_{I_k  }(x),\ x \in [0, 1],$ we define $\{L_{I }(x): x \geq 0  \}$ by gluing together the parts, more precisely by
\begin{eqnarray}\label{hnkhjs}
L_{I }(x) = L_{I_{\lfloor x \rfloor} }(x-\lfloor x\rfloor)   + \sum_{k=0}^{\lfloor x\rfloor-1} L_{I_k  }(1),\ \ \ \ \textrm{for all} \ x\geq 0.
\end{eqnarray}
Similarly, for $x< 0,$ we can define $L_{I }(x)=L_{I }(-x),$ since the function $\beta(x)=\alpha (-x)$
is defined on $[0, +\infty).$

\item   Let $(\phi(n))_{n \in \mathbb{N}}$ be a sequence of numbers satisfying $\phi(n)\rightarrow \infty$ as $n\rightarrow \infty.$ Assume that $\alpha(u)$ is  continuous in $[0,1].$ By an argument similar to
the proof of Theorem \ref{fhnk}, the sequence of processes
\begin{eqnarray}\label{fdvs}
\widehat{L}_{I } ^{(n)} (u)=   \sum_{k=1}^{\lfloor \phi(n) u \rfloor} \Big(\frac{1}{\phi(n)}\Big)^{1/\alpha (\frac{k}{\phi(n)})} X(k,n)\ , \ \ \  \ \ u \in [0,1],
\end{eqnarray}
tends in distribution to $L_I$ in $(D[0,1],d_S).$
Since $\alpha(u)$ is  continuous,  it is easy to see that $$\alpha \left(\frac{\lfloor \phi(n)u \rfloor}{\phi(n)} \right)\rightarrow\alpha (u)\ \ \ \ \textrm{as} \  n\rightarrow \infty.  $$  By the fact
that the summands of (\ref{fdvs}) verify
$$\Big(\frac{1}{\phi(n)}\Big)^{1/\alpha\big(\frac{\lfloor \phi(n)u \rfloor}{\phi(n)}\big)} X(\lfloor \phi(n)u \rfloor,n)\sim S_{\alpha(\frac{\lfloor \phi(n)u \rfloor}{\phi(n)}) }\Bigg( \Big(\frac{1}{\phi(n)}\Big)^{1/\alpha\big(\frac{\lfloor \phi(n)u \rfloor}{\phi(n)}\big)}, 0 , 0 \Bigg),$$
 equality (\ref{fdvs}) means that  the increment at the point $u$ of an  $\alpha(u)-$multistable process $L_{I}  (u)$ behaves locally like an $\alpha (u)-$stable random variable, but with the stability index $\alpha (u)$ varying with $u$.

\item
If $\alpha(u)\equiv\alpha$ for a constant $\alpha \in  (0,2]$, then $L_{I} (u)$ is just the usual
symmetric  $\alpha-$stable L\'{e}vy motion $L_{\alpha} (u)$. Hence, inequality (\ref{fdvs}) gives an equivalent definition
of the symmetric $\alpha-$stable L\'{e}vy motions: there is a sequence of independent and identically distributed (i.i.d.)
symmetric $\alpha-$stable random variables $(Y_{k})_{k \in \mathbb{N}}$ with an unit scale parameter such that
\begin{eqnarray}
L_{\alpha}^{(n)}  (u)  =   \sum_{k=1}^{\lfloor nu \rfloor} \frac{1}{n^{1/\alpha}} Y_{k} \ , \ \ \  \ \ u \in [0,1],
\end{eqnarray}
tends in distribution to $L_\alpha$ in $(D[0,1],d_S).$  This result is known as stable functional central limit theorem.

\item A slightly different method to construct
$L_{I}(u)$ can be stated as follows. Assume that
$\big( X( \frac{k}{ 2^n}) \big)_{n \in \mathbb{N}, \ k=1,...,2^n} $ is a family of independent and  symmetric $\alpha(\frac{k}{2^n})-$stable random variables with the unit scale parameter. Then it holds
\begin{eqnarray}\label{scv}
L_{I }(u)  =     \lim_ {n\rightarrow \infty }   \sum_{k=1}^{\lfloor 2^n u \rfloor} \Big(\frac{1}{2^n}\Big)^{1/\alpha(\frac{k}{2^n})} X\Big( \frac{k}{ 2^n}\Big) , \ \ \  \ \ \ \ \ \ \ u  \in [0,1],
\end{eqnarray}
where convergence is in  $(D[0,1],d_S).$   To highlight the differences between the two methods (\ref{multilevyI}) and (\ref{scv}), note that $X( \frac{k}{ 2^n})=X( \frac{2k}{ 2^{n+1}})$, while $X( k,  n )$ and $X( 2k, n+1)$ are two i.i.d. random variables.

\item  Inspecting the construction of field based MsLM  in Falconer and L\'{e}vy V\'{e}hel \cite{FL09},  it seems that the sequence of processes
\begin{eqnarray}
L^{(n)}_{F}(u)  {= } \sum_{k=1}^{\lfloor 2^n u \rfloor} \Big(\frac{1}{2^n}\Big)^{1/\alpha_n(u)} X(k,n) , \ \ \  \ \ \ \ \ \ \ u  \in [0,1],
\end{eqnarray}
tends in distribution to $L_F(u)$ in $(D[0,1],d_S)$. Unfortunately, it is not true in general.   We have the following  counter example.
\vspace{0.2cm}

\textbf{\emph{Example 1.}}   Consider the case of $\alpha_n(u)=\alpha(u)=\frac b 2 \mathbf{1}_{\{0\leq u\leq \frac b2 \} } + u \mathbf{1}_{\{ \frac b2 < u \leq 1 \}}.$ The characteristic function  of $L^{(n)}_{F}(u)$ is given by the following equality:  for any $\theta \in \mathbb{R},$
\begin{eqnarray}
 \mathbb{E} e^{i \theta L^{(n)}_{F}(u)}
 &{= }&  \prod_{k=1}^{\lfloor 2^n u \rfloor} \mathbb{E}\exp\left\{i  \theta \Big(\frac{1}{2^n}\Big)^{1/\alpha (u)} X(k,n) \right\} \nonumber\\
 &{= }& \exp\left\{-\sum_{k=1}^{  \lfloor2^n u\rfloor  } |\theta|^{  \alpha( \frac{k}{ 2^n})} \Big(\frac{1}{2^n} \Big)^{ \alpha(\frac{k}{ 2^n})/\alpha(u)}  \right\} ,\ \ \ \ \ u  \in [0,1].
\end{eqnarray}
Since,  for all $u \in (\frac b 2, 1]$ and $\theta\neq 0$,
\begin{eqnarray}
 \sum_{k=1}^{  \lfloor2^n u\rfloor  } |\theta|^{  \alpha(\frac{k}{ 2^n})} \Big(\frac{1}{2^n} \Big)^{ \alpha(\frac{k}{ 2^n})/\alpha(u)}  \geq  \sum_{k=1}^{  \lfloor2^n b/2\rfloor  } |\theta|^{b/2}\Big(\frac{1}{2^n} \Big)^{  b/2u } \rightarrow \infty, \ \ \ \ \ n\rightarrow \infty,
\end{eqnarray}
  we have $L^{(n)}_{F}(u)\rightarrow 0$ for all $u \in (\frac b 2, 1]$. Thus $L^{(n)}_{F}(u)$  does not tend  in distribution to $L_F(u)$ in  $(D[0,1],d_S).$

\end{enumerate}

\end{remark}

\section{Continuous Approximation of  MsLM}\label{contapp}
It is easy to see that when $\alpha(u)$ is a constant, then  the independent-increments MsLM reduce to $\alpha-$stable L\'{e}vy motions.  It is well known that $\alpha-$stable L\'{e}vy motions are stochastic H\"{o}lder continuous but not continuous.
We wonder if there exists a continuous approximation of   independent increments MsLM?  The answer is yes.

\subsection{A continuous stable process}
First, we shall construct a continuous stable process.   To this end, we shall make use of the following useful theorem.
\begin{theorem}\label{theo1} If the i.i.d.\ random variables $(Z_{jk})_{j,k}$ follow an $\alpha-$stable law, then it holds, for all $c> 1/\alpha,$
\[
\mathbb{P}\left( \bigcup_{i=1}^{\infty} \bigcap_{j\geq i}^{\infty} \max_{ k=0,...,2^{j}-1 } |Z_{jk}| \leq 2^{jc}\right)=1.
\]
\end{theorem}
\noindent\emph{Proof.} We only need to show that, for all $c> 1/\alpha$,
\[
\mathbb{P}\left( \bigcap_{i=1}^{\infty} \bigcup_{j\geq i}^{\infty} \max_{ k=0,...,2^{j}-1 } |Z_{jk}| > 2^{jc}\right)=0.
\]
By Borel-Cantelli Lemma, it is sufficient to prove, for all $c> 1/\alpha$,
\begin{equation}\label{ssd}
\sum_{j\geq 1} \mathbb{P}\left(  \max_{ k=0,...,2^{j}-1 } |Z_{jk}| > 2^{jc}\right) < \infty.
\end{equation}
To prove (\ref{ssd}), we need the following technical lemma  (cf.\ Property 1.2.15 of \mbox{Samorodnitsky} and Taqqu  \cite{ST94}   for details).
\begin{lemma} \label{lem1} Let $Z\sim S_\alpha(\sigma,\beta,\mu)$ with $0 <\alpha < 2.$ Then
\begin{displaymath}
\left\{ \begin{array}{ll}
\lim_{\lambda \rightarrow \infty} \lambda^{\alpha} \mathbb{P}( Z> \lambda) &=\  C_\alpha \frac{1+\beta}{2}\sigma^\alpha, \\
\\
\lim_{\lambda \rightarrow \infty} \lambda^{\alpha} \mathbb{P}( Z< - \lambda) & =\  C_\alpha \frac{1-\beta}{2}\sigma^\alpha.
\end{array} \right.
\end{displaymath}
\end{lemma}
Return to the proof of (\ref{ssd}). For all $c> 1/\alpha$ and all $j$ large enough, we have
\begin{eqnarray}
\mathbb{P}\left(  \max_{ k=0,...,2^{j}-1 } |Z_{jk}| > 2^{jc}\right) &=& 1-  \mathbb{P}\left(  |Z_{jk}| \leq 2^{jc}\ \textrm{for all}\ k=0,...,2^{j}-1 \right) \nonumber\\
 &=& 1 - \prod_{k=0}^{2^{j}-1} \mathbb{P}\left(  |Z_{jk}| \leq 2^{jc} \right).  \label{sns2}
\end{eqnarray}
Then, by equality (\ref{sns2}) and Lemma \ref{lem1}, we deduce
\begin{eqnarray*}
\mathbb{P}\left(  \max_{ k=0,...,2^{j}-1 } |Z_{jk}| > 2^{jc}\right)
 &=& 1 -   \left( 1+ O\Big(\frac { 1} {  2^{ j\alpha c} } \Big)  \right)^{2^{j}} \ \  \\
&=& O\Big(\frac{1}{ 2^{ j(\alpha c -1)}  }\Big),\ \ \ j\rightarrow \infty.
\end{eqnarray*}
Thus we obtain (\ref{ssd}) for all $c> 1/\alpha$.\qed

In the following theorem, we give a construction of continuous stable process. First, we recall the definition of the ``triangle'' function:
\begin{displaymath}
\varphi(t) = \left\{ \begin{array}{ll}
2t & \textrm{\ \ for $t \in [0, 1/2)$}\\
2-2t  & \textrm{\ \ for $t \in [1/2, 1]$}\\
0& \textrm{\ \ otherwise.}
\end{array} \right.
\end{displaymath}
Define $\varphi_{jk}(t)=\varphi(2^{j}t-k),$ for $j=0,1,...,$  and $k=0,...,2^{j}-1.$
\begin{theorem}\label{th2}
Assume the i.i.d.\ random variables $(Z_{jk})_{j,k}$ follow a symmetric $\alpha-$stable law with the unit scale parameter. Then, for all $d> 1/\alpha$, the process
$$X(t)  =   \sum_{j=0}^{\infty}\sum_{k=0}^{2^{j}-1}2^{-jd}Z_{jk}\varphi_{jk}(t),\ \ \ \ t \in [0,1],$$
 is a continuous and symmetric $\alpha$-stable process.  When $d=1/\alpha$, the process $X(t)$ is also a  symmetric, may not be continuous, $\alpha$-stable process
in $L^p(\Omega\times [0,1])$ for any $0<p<\alpha$.
\end{theorem}
\emph{Proof.} Set $X_{-1}\equiv 0$ and define the sequence of processes $(X_j)_{j \in \mathbb{N}}$ by:
\[
X_j(t)=X_{j-1}(t)+\sum_{k=0}^{2^{j}-1}2^{-jd}Z_{jk}\varphi_{jk}(t).
\]
First we show that the sequence of processes $(X_j)_{j \in \mathbb{N}}$ converges almost surely uniformly. Indeed, for all $t$,
\[
X_j(t)-X_{j-1}(t)=\sum_{k=0}^{2^{j}-1}2^{-jd}Z_{jk}\varphi_{jk}(t).
\]
Since the functions $(\varphi_{jk} )_{j,k}$ have disjoint supports and $|\varphi_{jk} |\leq 1$, it follows that
\[
|| X_j(t)-X_{j-1}(t) ||_{\infty}=2^{-jd} \max_{k=0,...,2^{j-1}} |Z_{jk}|.
\]
Theorem \ref{theo1} entails that $(X_j)_{j \in \mathbb{N}}$ converges almost surely in $C([0,1], || \cdot||_{\infty})$ to a continuous process $X$ for all $d> 1/\alpha.$  When $d=1/\alpha$,  we show that the sequence $(X_j)_{j \in \mathbb{N}}$ converges to a random variable $X$ in $L^p(\Omega\times [0,1])$ for any $0<p<\alpha$. Indeed, for any $0<p<\alpha$,
\begin{eqnarray}
\int_0^1  \mathbb{E}|X_j(t)-X_{j-1}(t)|^p dt &\leq& 2^{-jp/\alpha}\mathbb{E}|Z_{00}|^p \sum_{k=0}^{2^{j}-1}\int_0^1 \varphi_{jk}^p(t)dt  \nonumber \\
&\leq&  2^{-jp/\alpha}\mathbb{E}|Z_{00}|^p\int_0^1 \varphi_{00}^p(t)dt \nonumber \\
&=&  2^{-jp/\alpha}\mathbb{E}|Z_{00}|^p ,
\end{eqnarray}
this entitles convergence of $(X_j)_{j \in \mathbb{N}}$.

 Next, we prove that $X$ is a symmetric $\alpha-$stable process. By Theorem 3.1.2 of Samorodnitsky and Taqqu (1994),
we only need to check that all linear combinations
\[
\sum_{k=1}^d b_kX(t_k),\ \  d\geq 1,\ \  t_1,...,t_d \in [0,1]\  \textrm{and}\   b_1,...,b_d\  \textrm{real}
\]
are  symmetric $\alpha-$stable. We distinguish two cases as follows.  Define  $$D_n =\left\{ \frac{ k}{2^n}
 : 0 \leq  k \leq 2^n \right\}$$ and $D=\bigcup_{n=0,1,...} D_n$.\\
\textbf{i)}
If $t_k \in D, $ then all random variables $X(t_k), 1\leq k \leq d,$  are symmetric $\alpha-$stable. Thus all linear combinations $\sum_{k=1}^d b_kX(t_k)$ are  symmetric and $\alpha-$stable.\\
\textbf{ii)} For  $t_k \in [0,1], 1\leq k \leq d,$ we have $t_{kl} \in D$ such that $t_{kl} \rightarrow t_k, l\rightarrow \infty$. Since $X$ is continuous, we have $$\sum_{k=1}^d b_kX(t_k)= \lim_{j\rightarrow \infty} \sum_{k=1}^d b_kX(t_{kl}).$$
Its characteristic function has the following form:
\begin{eqnarray*}
\mathbb{E} \exp\Bigg\{i\theta \sum_{k=1}^d b_kX(t_k) \Bigg\} &=& \lim_{l\rightarrow \infty} \mathbb{E} \exp\Bigg\{i\theta \sum_{k=1}^d b_kX(t_{kl}) \Bigg\}.
\end{eqnarray*}
It is easy to see that  the scale parameter of $\sum_{k=1}^d b_kX(t_{kl})$ is
$$  \sigma_{l}(\alpha)=\left( \sum_{j=0}^{\infty}\sum_{i=0}^{2^{j}-1} \Big( \sum_{k=1}^d |b_k| 2^{-jd} \varphi(2^jt_{kl}-i)\Big)^\alpha \right)^{1/\alpha} .$$
Since at most one summand of the sum $\sum_{i=0}^{2^{j}-1}   2^{-jd} \varphi(2^jt-i)$ is non-zero and  $$ \sum_{i=0}^{2^{j}-1} \Big(\sum_{k=1}^d |b_k| 2^{-jd} \varphi(2^jt_{kl} -i)\Big)^\alpha  \leq d\, b^\alpha   2^{-j\alpha d},  $$
 where $b=\max\{|b_k|, 1\leq k \leq d \},$
then $\sigma(\alpha) = \lim_{l\rightarrow \infty}\sigma_l(\alpha)$ exists for $d\geq 1/\alpha$ and
\begin{eqnarray*}
\mathbb{E} \exp\Bigg\{i\theta \sum_{k=1}^d b_kX(t_k) \Bigg\}
  &=& \lim_{j\rightarrow \infty} \exp\Big\{ - \sigma_l(\alpha)^{\alpha} |\theta|^{\alpha} \Big\} \ \  \\
  &=& \exp\Big\{ - \sigma(\alpha)^{\alpha} |\theta|^{\alpha} \Big\}.
\end{eqnarray*}
This implies that all linear combinations $\sum_{k=1}^d b_kX(t_k)$ are symmetric $\alpha-$stable random variables.   This completes the proof. \qed

One deduces the scale parameter $\sigma(t)$ of the process $X(t)$ is given as follows
\[
\sigma^\alpha(t)=  \sum_{j=0}^{\infty}\sum_{k=0}^{2^{j}-1} \Big( 2^{-jd} \varphi(2^jt -k)\Big)^\alpha.
\]
By noting that at most one $\varphi(2^jt -k)$ is non-zero for all $j$, we have the following estimation of the scale parameter
\[
 \varphi^{1/\alpha}( t ) \ \leq \ \sigma (t) \ \leq \ \Big(\frac{1}{1-2^{-\alpha d}} \Big)^{1/\alpha},\ \ \ \ t \in [0, 1].
\]
It is worth  noting that when $t\neq0,1,$  we have $\sigma (t)>0.$ This observation will be useful to establish
 continuous approximations of MsLM in the next subsection.

\subsection{Continuous approximations of MsLM}
In Theorem \ref{fhnk}, we establish discrete approximations of  the independent-increments MsLM.
In this subsection, we shall give continuous approximations of the  independent-increments MsLM.  It is worth to noting that one
cannot make use of the method of Theorem \ref{th2} to establish continuous approximations
of MsLM in general, since a sum  of two stable random variables
with different stability indices is not a stable random variable.  To obtaining continuous
approximations of  the independent-increments MsLM, our main method is to replace the summands in  (\ref{multilevyI})
 by a sequence of independent and continuous stable processes starting at $0$, for instance the stable processes established in Theorem \ref{th2}.
\begin{theorem}\label{cosa} Let $\alpha(u)$ be a continuous function ranging in $[a,b] \subset (0,2]$.
Assume that $\left(X_{\alpha(\frac{k}{2^n})}( t) \right)_{ n \in \mathbb{N}, \, k=0,...,2^n-1}$
 is a family of independent and  continuous   $\alpha(\frac{k}{2^n})-$stable random processes. Assume  $X_{\alpha(\frac{k}{2^n})} ( 0)=0$ and $\sigma_{\alpha(\frac{k}{2^n})}(t)>0$ for all $t\in (0, 1]$ and all $n \in \mathbb{N}, \, k=0,...,2^n-1$, where $\sigma_{\alpha(\frac{k}{2^n})}(t)$ is the scale parameter of  $X_{\alpha(\frac{k}{2^n})} (t)$.   Define
\begin{eqnarray}\label{clevyc}
S_n(u)  &  = &  \Big(\frac{1}{2^n}\Big)^{\alpha(\frac{\lfloor 2^n u \rfloor}{2^n})} \frac{1}{\sigma_{\alpha(\frac{\lfloor 2^n u \rfloor}{2^n})}(\frac{1}{2^n})}  X_{\alpha(\frac{\lfloor 2^n u \rfloor}{2^n})} \Bigg(u- \frac{\lfloor 2^n u \rfloor}{2^n}   \Bigg) \nonumber\\
&& \ \ \ \  +  \sum_{k=0}^{\lfloor 2^n u \rfloor -1} \Big(\frac{1}{2^n}\Big)^{\alpha(\frac{k}{2^n})} \frac{1}{\sigma_{\alpha(\frac{k}{2^n})}(\frac{1}{2^n})}  X_{ \alpha(\frac{k}{2^n})}  \Bigg( \frac{1}{2^n}   \Bigg), \ \ \  \ \ u  \in [0,1].
\end{eqnarray}
Then $(S_n)_{n \in \mathbb{N}}$ is a sequence of continuous processes and  the process $S_n(u),  u  \in [0,1],$
tends in distribution to $L_I(u)$ in $(D[0,1],d_S).$
\end{theorem}

By the definition of $S_n(u)$ in (\ref{clevyc}), it seems that the process $S_n(u)$ restores more and more details of $L_{I }(u)$ when $n$ is increasing.

It is worth noting that when $\alpha(u)\equiv \alpha$ for a constant $\alpha \in (0, 2]$, Theorem \ref{cosa} gives continuous approximations to the usual symmetric $\alpha-$stable L\'{e}vy motion $L_\alpha(u).$

\noindent\emph{Proof.} It is easy to see that  the first item in the right hand side of (\ref{clevyc}) converges to zero in distribution as
$n\rightarrow \infty$, i.e.,
\begin{eqnarray}
\lim_{n\rightarrow \infty}\left( \Big(\frac{1}{2^n}\Big)^{\alpha(\frac{\lfloor 2^n u \rfloor}{2^n})} \frac{1}{\sigma_{\alpha(\frac{\lfloor 2^n u \rfloor}{2^n})}(\frac{1}{2^n})}  X_{\alpha(\frac{\lfloor 2^n u \rfloor}{2^n})} \Bigg(u- \frac{\lfloor 2^n u \rfloor}{2^n}   \Bigg) \right) = 0
\end{eqnarray}
in distribution. Notice that the summands
\begin{eqnarray}
 \frac{1}{\sigma_{\alpha(\frac{k}{2^n})}(\frac{1}{2^n})}  X_{ \alpha(\frac{k}{2^n})}  \Bigg( \frac{1}{2^n}   \Bigg)
\end{eqnarray}
in the right hand side of (\ref{clevyc}) are independent $\alpha(\frac{k}{2^n})-$stable random variables with the unit scale parameter. Using Theorem \ref{fhnk}, we find that the process $S_n(u),  u  \in [0,1],$
tends in distribution to $L_I(u)$ in $(D[0,1],d_S).$ \qed

\section{Integrals of Multistable L\'{e}vy Measure}\label{endsection}
Let $\alpha=\alpha(u), u  \in [0,1],$ be a c\`{a}dl\`{a}g function
ranging in $[a,b] \subset (0,2].$  Denote by
\[
\mathcal{L}_{\alpha }[0,1] =\Big\{ f:  f \textrm{ is measurable with }  ||f||_{\alpha  } < \infty \Big\},
\]
where
\[
||f||_{\alpha  }
:= \inf \left\{  \lambda > 0: \int_{0}^{1}\Big|\frac{f(x)}{\lambda} \Big|^{ \alpha(x)} dx =1  \right\} \ \ \ \ \ \textrm{and} \ \ \ \ \  ||0||_{\alpha}=0.
\]
Note that $||\cdot||_{\alpha }$  is a quasinorm; see   Falconer and  Liu \cite{FL12} and Ayache \cite{AA13}.
Using the Kolmogorov consistency conditions and the L\'{e}vy continuity theorem, Falconer and  Liu \cite{FL12} (see also Falconer \cite{F09})  proved that the characteristic function, for all $(\theta_1,...,\theta_d)\in \mathbb{R}^d,$
\begin{eqnarray}\label{fli1}
\mathbb{E}  \exp\left\{ i \Bigg(\sum_{j=1}^d \theta_j I(f_j)   \Bigg)  \right\}  =\exp \left\{- \int \Big| \sum_{j=1}^d \theta_j f_j(x) \Big|^{\alpha(x)} dx \right\} \ \ \
\end{eqnarray}
well defines a consistent probability distribution of the random vector $(I(f_1), I(f_2),...,I(f_d))\in \mathbb{R}^d$ on the functions $f_j \in \mathcal{L}_{\alpha }[0,1],$  where $I(f)=\int f(x) M_\alpha(d x).$   They called $M_\alpha$ the multistable L\'{e}vy measure and
$I(f)=\int f(x) M_\alpha(d x)$ the integral with respect to $M_\alpha.$
Moreover, they also showed that   the integrals of functions with disjoint supports are independent.
In particular, it holds $$L_I(u)=\int   \mathbf{1}_{[0,\ u]}(x)M_\alpha ( dx),\ \ \ \ \   u  \in [0,1].$$

In the following theorem, we give an alternative definition of the integrals based on
the weighted sums of independent random variables.
\begin{theorem} Let  $\big( X(k,n) \big)_{n \in \mathbb{N}, \ k=1,...,2^n} $ be defined by Theorem \ref{fhnk}.   Then, for any $f \in \mathcal{L}_{\alpha }[0,1],$  it holds
\begin{eqnarray}
\int_0^1 f(x) M_\alpha(d x) =  \lim_ {n\rightarrow \infty }   \sum_{k=1}^{ 2^n  } \Big(\frac{1}{2^n}\Big)^{1/\alpha(\frac{k}{2^n})} f\Big(\frac{k}{2^n}\Big) X(k,n)
\end{eqnarray}
in distribution.
\end{theorem}
\emph{Proof.} Denote by
\[
 S(k,n)=\Big(\frac{1}{2^n}\Big)^{1/\alpha(\frac{k}{2^n})} f\Big(\frac{k}{2^n}\Big) X(k,n) \ \ \ \  and \ \ \ \ X_n = \sum_{k=1}^{  2^n  } S(k,n) .
\]
It is easy to see that, for any $\theta \in \mathbb{R},$
\begin{eqnarray*}
\mathbb{E} e^{i\, \theta X_n  } &=&   \prod_{k=1}^{ 2^n  }  \mathbb{E} e^{i \theta S(k,n) }
 =  \exp\Bigg\{-\Big| \theta f\Big(\frac{k}{2^n}\Big)\Big|^{\alpha(k/2^n)} \frac{1}{2^n} \Bigg\} .
\end{eqnarray*}
Hence, we have
\begin{eqnarray*}
\lim_{n\rightarrow \infty}\mathbb{E} e^{i\, \theta X_n  } &=&   \exp \left\{- \int_0^1 \Big|  \theta  f (x) \Big|^{\alpha(x)} dx \right\},
\end{eqnarray*}
which means $ \lim_ {n\rightarrow \infty  }  X_n =\int_0^1 f(x) M_\alpha(d x)$ in distribution by the definition  of the multistable integrals with respect to the multistable L\'{e}vy measure $M_\alpha.$  \qed

The following theorem relates the convergence of a sequence of $\alpha(u)-$multistable integrals
to the convergence of the sequence of integrands.

\begin{theorem}\label{scxfd} Assume $X_j=\int_0^1 f_j(x) M_\alpha(d x)$ and $X=\int_0^1 f(x) M_\alpha(d x),$
for $f_j, j=1,2,..., f \in \mathcal{L}_{\alpha  }[0,1].$ Then
\[
\lim_{j\rightarrow \infty}  X_j = X
\]
in probability, or
\[
\lim_{j\rightarrow \infty} ( X_j -X)= 0
\]
in distribution, if and only if
\[
\lim_{j\rightarrow \infty} || f_j   - f  ||_{\alpha  }   = 0.
\]
\end{theorem}
\emph{Proof.} The convergence $\lim_{j\rightarrow \infty} X_j = X
$ in probability   is equivalent to $\lim_{j\rightarrow \infty} (X_j -X)= 0
$ in probability   and hence to the convergence in distribution to zero of the sequence $(X_j -X)_{j=1,2,...}.$
If $X_j-X$ convergence in distribution to $0,$ then, for any $\theta \in \mathbb{R},$
\begin{eqnarray}
1= \lim_{j\rightarrow \infty}\mathbb{E} e^{i\, \theta (X_j -X)  }  =  \lim_{j\rightarrow \infty} \exp \left\{- \int_0^1 \Big|  \theta \Big(f_j (x) - f (x)\Big) \Big|^{\alpha(x)} dx \right\},
\end{eqnarray}
which is equivalent to, for any $\lambda>0,$
\[
\lim_{j\rightarrow \infty} \int_0^1 \Big| \frac{ f_j (x) - f (x)}{\lambda }   \Big|^{\alpha(x)} dx = 0.
\]
This equality means $\lim_{j\rightarrow \infty} || f_j   - f  ||_{\alpha  }   = 0.$
\qed

The last theorem shows that convergence in probability of multistable integrals coincides with convergence in quasinorm $||\cdot||_{\alpha  }$.

The convergence $\lim_{j\rightarrow \infty} X_j = X
$  almost surely  implies  the convergence $\lim_{j\rightarrow \infty} X_j = X
$ in probability.  Thus the following corollary is obvious.
\begin{corollary}
Assume that $X_j, j=1,2,..$  and $X$ are defined by Theorem \ref{scxfd}.
If
\[
\lim_{j\rightarrow \infty}  X_j  =X
\]
  almost surely, then
\[
\lim_{j\rightarrow \infty} || f_j   - f  ||_{\alpha  }   = 0.
\]
\end{corollary}

\subsection{Independence}
Independence of two multistable integrals imposes a stronger restriction on the  integrands: they
must almost surely have disjoint supports with respect to Lebesgue measure $\mathcal{L}$. Indeed,
\begin{theorem}\label{THinde} Let  $X_1= \int_0^1 f_1(x)M_\alpha(dx)$ and $X_2= \int_0^1 f_2(x)M_\alpha(dx)$
be two  multistable integrals, where  $f_j\in \mathcal{L}_{\alpha  }[0,1], j=1,2.$ Assume either
 $$[a,b] \subset (0, 2)$$ or
 \begin{eqnarray} \label{tonghao}
f_1(x)f_2(x) \geq 0\ \ \ \ \  \mathcal{L}-a.s.\textrm{ on }[0, \, 1].
\end{eqnarray}
Then $X_1$ and $X_2$ are independent if and only if
\begin{eqnarray}\label{dssdf}
f_1(x)f_2(x) \equiv 0\ \ \ \ \  \mathcal{L}-a.s.\textrm{ on }[0, \, 1].
\end{eqnarray}
\end{theorem}
\emph{Proof.}
Two multistable integrals $X_1$ and $X_2$ are independent if and only if, for any $(\theta_1, \theta_2) \in \mathbb{R}^2,$
\begin{eqnarray}\label{dfds}
\mathbb{E} \exp\Big\{i(\theta_1 X_1+ \theta_2 X_2 ) \Big\}
  &=& \mathbb{E} \exp\Big\{i\theta_1 X_1\Big\} \  \mathbb{E} \exp\Big\{i\theta_2 X_2\Big\}.
\end{eqnarray}
Notice that
\begin{eqnarray*}
\mathbb{E} \exp\Big\{i(\theta_1 X_1+ \theta_2 X_2 ) \Big\}  = \exp \left\{- \int_0^1 \Big| \sum_{j=1}^2 \theta_j  f_j(x)  \Big|^{\alpha(x)} dx \right\}, \ \ \
\end{eqnarray*}
and that
\begin{eqnarray*}
\mathbb{E} \exp\Big\{i\theta_1 X_1\Big\} \  \mathbb{E} \exp\Big\{i\theta_2 X_2\Big\}  =\exp \left\{- \sum_{j=1}^2\int_0^1 \Big|  \theta_j  f_j(x)  \Big|^{\alpha(x)} dx \right\}.
\end{eqnarray*}
Equating the moduli of (\ref{dfds}) gives
\begin{eqnarray}\label{sddvkf}
 \int_0^1 \Big| \sum_{j=1}^2 \theta_j  f_j(x)  \Big|^{\alpha(x)} dx= \sum_{j=1}^2\int_0^1 \Big|  \theta_j f_j(x) \Big|^{\alpha(x)} dx  .
\end{eqnarray}
Notice that (\ref{sddvkf}) implies that
\begin{eqnarray}
 \int_0^1 \Big|   f_1(x) - f_2(x)  \Big|^{\alpha(x)} dx &=&  \int_0^1 \Big|  f_1(x)   \Big|^{\alpha(x)} dx +   \int_0^1 \Big|  f_2(x)  \Big|^{\alpha(x)} dx.\label{sddhnlg2}\\
 &=&\int_0^1 \Big|  f_1(x) + f_2(x)  \Big|^{\alpha(x)} dx   \label{sddhnlg1}
\end{eqnarray}

Assume $[a,b] \subset (0, 2).$ We argue as Lemma 2.7.14 of Samorodnitsky and Taqqu \cite{ST94}.  When $\alpha  \in (0, 2),$  the function $r_\alpha(u)=u^{\alpha/2}, u\geq 0,$ is strictly concave. Therefore, for fix $x \in [0, 1],$
\begin{eqnarray}
&& |f_1(x)+f_2(x)|^{\alpha(x)}+ |f_1(x)-f_2(x)|^{\alpha(x)} \nonumber \\
&& \ \ \  \ \ \  \ \ \ \ \ \ = 2 \, \frac{r_{\alpha(x)}(|f_1(x)+f_2(x)|^2) + r_{\alpha(x)}(|f_1(x)-f_2(x)|^2)}{2} \nonumber\\
&& \ \ \  \ \ \  \ \ \ \ \ \ \leq 2 \, r_{\alpha(x)} \Big(  \frac{|f_1(x)+f_2(x)|^2 +  |f_1(x)-f_2(x)|^2}{2} \Big) \nonumber\\
&& \ \ \  \ \ \  \ \ \ \ \ \ = 2 \, r_{\alpha(x)} \big( f_1(x)^2    +  f_2(x)^2  \big) \nonumber\\
&& \ \ \  \ \ \  \ \ \ \ \ \ \leq  2 \, \big( |f_1(x)|^{\alpha(x)}+ |f_2(x)|^{\alpha(x)} \big) \label{fdd1ss}
\end{eqnarray}
with equality in the preceding relations equivalent $f_1(x)f_2(x)=0.$  Inequalities (\ref{sddhnlg2}) and (\ref{sddhnlg1}) imply that
\begin{eqnarray}
&& \int_0^1 \Big|   f_1(x) - f_2(x)  \Big|^{\alpha(x)} dx +   \int_0^1 \Big|  f_1(x) + f_2(x)  \Big|^{\alpha(x)} dx  \ \ \ \ \  \ \ \ \ \ \ \   \nonumber \\
 & & \ \ \ \ \  \ \ \ \ \ \ \ \ \ \ \ \  \ \ \ \ \ \ \    = 2 \Bigg( \int_0^1 \Big|  f_1(x)   \Big|^{\alpha(x)} dx +   \int_0^1 \Big|  f_2(x)  \Big|^{\alpha(x)} dx \Bigg).
\end{eqnarray}
Now (\ref{fdd1ss}) implies that the left-hand side of the last equality is always less than or equal to the right-hand side of the last inequality and, if
they are equal, then necessarily  (\ref{dssdf}) holds.

Assume (\ref{tonghao}). Then it holds $\big|   f_1(x) - f_2(x)  \big| \leq \big|  f_1(x) + f_2(x)  \big| \ \mathcal{L}-a.s.\textrm{ on }[0, \, 1].$ When $\alpha  \in (0, 2],$  the function $r_\alpha(u)=u^{\alpha/2}$ is increasing in $u \in [0, \infty).$ Hence (\ref{sddhnlg1})   holds if and only if (\ref{dssdf}) holds.

This proves that (\ref{dssdf}) is a necessary condition for the independence of $X_1$ and $X_2.$
It is also sufficient because if (\ref{dssdf}) holds, then  (\ref{sddvkf}) also holds.  \qed

The preceding result is very useful and will often be used in the sequel.
\begin{theorem}\label{Tdsde} Assume $f_j\in \mathcal{L}_{\alpha  }[0,1],$ $j=1,...,d.$
Assume either $[a,b] \subset (0, 2)$  or $f_i(x)f_k(x) \geq 0$   $\ \mathcal{L}-a.s.\textrm{ on }[0, \, 1]$ for any subset $\{i, k\}$ of $\{1,2,..., d\}.$
The multistable integrals  $X_j= \int_0^1 f_j(x)M_\alpha(dx),$  $j=1,...,d,$
are independent if and only if they are pairwise independent, i.e., if and only if
 \begin{eqnarray}\label{gbjs}
f_i(x)f_k(x)\equiv 0  \ \ \  \mathcal{L}-a.s.\textrm{ on }[0, \, 1]
\end{eqnarray}
for any subset $\{i, k\}$ of $\{1,2,..., d\}.$
\end{theorem}
\emph{Proof.} Independence clearly implies pairwise independence. By Theorem \ref{THinde}, pairwise independence  implies (\ref{gbjs}).
If (\ref{gbjs}) holds,  then it holds, for any $(\theta_1, ..., \theta_d) \in \mathbb{R}^d,$
\begin{eqnarray}
 \int_0^1 \Big| \sum_{j=1}^d \theta_j  f_j(x)  \Big|^{\alpha(x)} dx= \sum_{j=1}^d \int_0^1 \Big|  \theta_j f_j(x) \Big|^{\alpha(x)} dx  .
\end{eqnarray}
Thus  the joint characteristic  function of $X_1,..., X_d$ factorizes
\begin{eqnarray*}
\mathbb{E} \exp\Big\{i \sum_{j=1}^d \theta_j X_j  \Big\}
  &=& \exp \left\{- \int_0^1 \Big| \sum_{j=1}^d \theta_j  f_j(x)  \Big|^{\alpha(x)} dx \right\} \\
  &=& \exp \left\{- \sum_{j=1}^d\int_0^1 \Big| \theta_j  f_j(x)  \Big|^{\alpha(x)} dx \right\} \\
  &=& \prod_{j=1}^d \mathbb{E} \exp\Big\{i \theta_j X_j  \Big\}.
\end{eqnarray*}
This proves that $X_1,..., X_d$ are independent. \qed

\subsection{Stochastic H\"{o}lder continuity}
We call a random process $X(u),\, u \in  I,$ is  \emph{stochastic H\"{o}lder continuous} of exponent $\beta \in (0, 1],$ if it holds
\[
\limsup_{ u,r \in I,\  |u-r| \rightarrow   0 }   \mathbb{P}\big(|X(u)-X(r)|\geq C |u-r|^\beta \big)=0
\]
for a positive constant $C.$  It is obvious  that
if $X(u)$ is  stochastic  H\"{o}lder continuous of exponent $\beta_1 \in (0, 1],$ then $X(u)$ is
stochastic  H\"{o}lder continuous  of exponent $\beta_2 \in (0, \beta_1].$

\vspace{0.2cm}
\noindent\textbf{Example 2.}  \emph{Assume that a random process $X(u),\, u \in  I,$ satisfies the following condition: there exist  three
strictly positive constants $\gamma, c, \rho$ such that
\[
\mathbb{E} |X(u)-X(r)|^\gamma    \leq c \, |u-r|^\rho,\ \ \ \ \ \ \  u,r \in I.
\]
Then $X(u), u \in  I,$ is  stochastic H\"{o}lder continuous  of exponent $\beta \in (0, \min\{1, \rho/\gamma\} ).$ Indeed, it is easy to see that for all $u, r \in I,$
\begin{eqnarray*}
\mathbb{P}\Big( |X(u)-X(r)|\geq C |u-r|^\beta \Big) &\leq& \frac{\mathbb{E} |X(u)-X(r)|^\gamma   }{C^\gamma |u-r|^{\beta \gamma}} \\
&\leq& \frac{c}{C^\gamma}  |u-r|^{\rho-\beta \gamma},
\end{eqnarray*}
which implies our claim.  }

\vspace{0.2cm}

The following theorem gives a sufficient condition such that the integrals with respect to multistable L\'{e}vy measure $M_\alpha$ are stochastic H\"{o}lder continuous.
\begin{theorem}\label{endthm} Assume that  $X(t)= \int_0^1 f(t,x)M_\alpha(dx)$ is a  multistable integral,
where $f(t,x)$ is jointly measurable and $f(t,x)  \in \mathcal{L}_{\alpha }[0,1]$ for all $t\in I.$
If there exist two constants $\eta > 0 $ and $C>0$ such that
\begin{eqnarray}\label{fine37}
  \int_0^1 \Big|
 f(t,s)-f(v,s)   \Big|^{\alpha(s)} ds
&\leq&    C \, \Big| t - v \Big|^\eta, \ \ \ \  \ \ t, v \in I.
\end{eqnarray}
Then it holds
\begin{eqnarray}\label{finein34}
 \mathbb{P}( |X (t)-X (v)| \geq | t - v  |^\beta) \leq C_{a,b}\ | t - v  |^{ \eta- b\beta},\ \ \ \  \ \ t, v \in I,
\end{eqnarray}
where $C_{a,b}$ is a constant depending on $a, b$ and $C.$ In particular, it implies that
$X(t)$  is stochastic H\"{o}lder continuous of exponent $\beta \in (0, \min\{1, \eta/b\})$.

\end{theorem}

\noindent\emph{Proof.}  By the Billingsley inequality (cf.  p.\ 47 of \cite{B68}) and (\ref{fine37}), it is easy to see that, for all  $t, v \in  I   $   and all $x > 0,$
\begin{eqnarray*}
 \mathbb{P}( |X (t)-X (v)| \geq x)
 &\leq&  \frac x2 \int_{-2/x}^{\,2/x} \left( 1 -  \mathbb{E}e^ {i \theta \big(X (t)-X (v) \big)}\right)\, d\theta\\
&=&  \frac x2 \int_{-2/x}^{\,2/x} \left( 1 -    \exp \left\{- \int_0^1 \Big|  \theta  \Big( f(t,z)-f(v,z) \Big)\Big|^{\alpha(z)} dz \right\} \right)\, d\theta\\
&\leq&  \frac x2 \int_{-2/x}^{\,2/x} \int_0^1 \Big|  \theta  \Big( f(t,z)-f(v,z) \Big)\Big|^{\alpha(z)} dz  \, d\theta\\
&\leq&  \frac x2 \Bigg[\int_{|\theta|< 1}\Big|  \theta\Big|^{a}      \, d\theta  +\int_{1\leq |\theta| \leq 2/x} \Big|  \theta \Big|^{b}  d\theta \Bigg] \  C \Big| t - v \Big|^\eta\\
&\leq&   C \left( \frac{x\,  }{a+1}  +   \frac{2^{b+1}   } {b+1 } \frac{1} {x^b } \right) \Big| t - v \Big|^\eta.
\end{eqnarray*}
Taking $x= | t - v  |^\beta,$ we obtain (\ref{finein34}).
This implies that $X(t)$  is stochastic H\"{o}lder continuous of exponent $\beta \in (0, \min\{1, \eta/b\})$.
 \qed

As an example to illustrate Theorem \ref{endthm}, consider the weighted MsLM introduced by  \mbox{Falconer}  and  Liu \cite{FL12}. The following theorem shows that  the weighted  MsLM are  H\"{o}lder continuous of exponent $\beta \in (0, \min\{1, 1/b\})$.
\begin{theorem} \label{lemma2}  Let
 \[
Y(t)=\int_0^1 w(x)\mathbf{1}_{[0,\ t]}(x)M_\alpha ( dx),\ \ \ t  \in [0,1],
\]
 be a weighted multistable L\'{e}vy motion, where the function $ w(x), x \in [0, 1],$ is c\`{a}dl\`{a}g.
Then   $Y(t)$ is stochastic H\"{o}lder continuous of exponent $\beta \in (0, \min\{1, 1/b\}).$ Moreover,
 it holds
\begin{eqnarray} \label{fskfdf}
 \mathbb{P}( |
 Y (t)-Y (v)| \geq | t - v  |^\beta) \leq C_{a,b}\ | t - v  |^{ 1- b\beta},\ \ \ \  \ \ t, v  \in [0,1],
\end{eqnarray}
where $C_{a,b}$ is a constant depending on $a, b, \alpha(\cdot)$ and $w(\cdot).$
In particular,  it implies that $L_{I } (u), u  \in [0,1],$ is stochastic H\"{o}lder continuous of exponent $\beta \in (0, \min\{1, 1/b\})$.
\end{theorem}

\noindent\emph{Proof.} Set $f(t,x)=w(x)\mathbf{1}_{[0, t]}(x), \ t,x \in [0, 1].$ It is easy to see that, for all $v, t \in [0,1]$ such that $v \leq t,$
\begin{eqnarray*}
 \int_0^1 \Big|
 f(t,s)-f(v,s)   \Big|^{\alpha(s)} ds&\leq&  \int_0^1 \Big|
w(s)\mathbf{1}_{[v,\  t]}(s) \Big|^{\alpha(s)} ds\\
&\leq& C_\omega \,  \int_0^1
 \mathbf{1}_{[v,\  t]}(s)   \, ds   \\
&\leq&  C_\omega \,  (t-v),
\end{eqnarray*}
where $C_\omega = \sup_{z \in [0, 1]}|w(z)|^{\alpha(z)}.$ By Theorem \ref{endthm}, we get (\ref{fskfdf}). This completes the proof of Theorem \ref{lemma2}. \qed

\subsection{Strongly localisability}

When the function $\alpha(x) \in [a, b], x \in [0, 1], $ is continuous,  some sufficient conditions such that the multistable integrals
 are localisable (or strongly localisable) has been obtained by  Falconer  and  Liu.
In the following theorem, we give some new conditions  such that  localisability  can be strengthened to strongly localisability.
\begin{theorem}\label{scrend}  Assume that $f(t,x)$ and $h(t, x)$  are jointly measurable; and that $f(t,x), h(t, x) \in \mathcal{L}_{\alpha }[0,1]$ for any $t\in [0, 1]$.  Assume that  $X(t)= \int_0^1 f(t,x)M_\alpha(dx)$ and $X_x'(t)=\int_0^1 h(t,x) M_\alpha(d x)$ are two multistable integrals and have versions in $D[0, 1]$.
Suppose that $X(t)$ is $1/\alpha(x)-$localisable at $x$ with local form $X_x'(t).$
If there exist two constants $\eta > 1 $ and $C>0$ such that
\begin{eqnarray}\label{vbsds}
  \int_0^1 \Bigg| \frac{
 f(x+rt,s)-f(x+rv,s)  }{ r^{1/\alpha(x)} } \Bigg|^{\alpha(s)} ds
&\leq&    C \Big| t - v \Big|^\eta,\ \ \ \  \ \ \ t, v \in [0, 1],
\end{eqnarray}
for all sufficiently small $r>0,$ then $X(t)$  is strongly localisable at all $x \in [0, 1].$
Moreover, if $X(t)$ has independent increments and $(\ref{vbsds})$ holds for a constant $\eta > 1/2$, then  the claim holds also.
\end{theorem}

Notice that  condition (\ref{vbsds}) is slightly more general than the condition of Falconer  and  Liu  (cf.\ Theorem 3.2 of \cite{FL12}): there exist two constants $\eta > 1/a $ and $C>0$ such that
\begin{eqnarray}\label{vbsds1}
 \Bigg|\Bigg| \frac{
 f(x+rt,\cdot)-f(x+rv,\cdot)  }{ r^{1/\alpha(x)} } \Bigg|\Bigg|_{\alpha }
&\leq&    C \Big| t - v \Big|^\eta,\ \ \ \ \ \ \  t, v \in [0, 1],
\end{eqnarray}
for all sufficiently small $r>0.$

\noindent\emph{Proof.} For any $x \in [0, \, 1),$ define $$X_r(u)=\frac{X(x+ru )-X(x )}{ r^{1/\alpha(x)} } \,,\ \ \ \  \ \  r, u \in (0, \, 1].$$
By Theorem 15.6 of \mbox{Billingsley \cite{B68}}, it suffices to show that, for some $\beta> 1$ and $\tau\geq 0$,
\begin{eqnarray}\label{ineq38}
\mathbb{P}\Big( \Big| X_r(u)-X_r(u_1)\Big|\geq \lambda,\ \Big|X_r(u_2)-X_r(u)\Big|\geq \lambda \Big) \leq \frac{C}{\lambda^{ \tau } } \Big[  u_2 - u_1 \Big]^{ \beta}
\end{eqnarray}
for $u_1 \leq u \leq u_2, \lambda> 0$ and $ r \in (0, 1]$, where
$C$ is a positive constant.
Since $X_r(u)-X_r(u_1)$ and $X_r(u_2)-X_r(u)$
 are symmetric,
it follows that
\begin{eqnarray*}
&& \mathbb{P}\Bigg( \Big| X_r(u)-X_r(u_1)\Big|\geq \lambda,\ \Big|X_r(u_2)-X_r(u)\Big|\geq \lambda \Bigg) \\
&\leq&4 \, \mathbb{P}\Bigg( X_r(u)-X_r(u_1) + \Big( X_r(u_2)-X_r(u) \Big)\geq 2\lambda   \Bigg)  \\
&=&4 \, \mathbb{P}\Big( X_r(u_2)-X_r(u_1) \geq 2\lambda   \Big).
\end{eqnarray*}
By the Billingsley inequality (cf.  p.\ 47 of \cite{B68}) and  (\ref{vbsds}), we have
\begin{eqnarray}
\mathbb{P}\Big( X_r(u_2)-X_r(u_1) \geq 2\lambda   \Big) &\leq&  \lambda  \int_{-1/\lambda}^{1/\lambda } \Bigg( 1- \mathbb{E}e^{i\theta \big(X_r(u_2)-X_r(u_1)\big)}\Bigg) d\theta \nonumber \\
&=&  \lambda  \int_{-1/\lambda}^{1/\lambda } \Bigg( 1-  e^{- \int_0^1 \big|\theta  \frac{
 f(x+ru_2,s)-f(x+ru_1,s)  }{ r^{1/\alpha(x)} }  \big|^{\alpha(s)} ds   }\Bigg) d\theta \nonumber  \\
&\leq&  \lambda  \int_{-1/\lambda}^{1/\lambda } \int_0^1 \Bigg|\theta  \frac{
 f(x+ru_2,s)-f(x+ru_1,s)  }{ r^{1/\alpha(x)} } \Bigg|^{\alpha(s)} ds \, d\theta \nonumber  \\
&=& \lambda  \int_{-1/\lambda}^{1/\lambda } |\theta|^{\mu} \int_0^1 \Bigg| \frac{
 f(x+ru_2,s)-f(x+ru_1,s)  }{ r^{1/\alpha(x)} } \Bigg|^{\alpha(s)} ds \, d\theta \nonumber \\
&\leq&\frac{C_1}{\lambda^{\gamma} } \Big[ u_2   - u_1  \Big]^\eta, \label{sfnmbs}
\end{eqnarray}
where  $\mu=a\mathbf{1}_{[1,\, \infty)} (\theta)+ b\mathbf{1}_{(0,\, 1 )}(\theta),    \gamma = a\mathbf{1}_{[1,\, \infty)} (\lambda)+ b\mathbf{1}_{(0,\, 1 )}(\lambda)$ and $C_1$ is  a positive constant  depending only on $a, b$ and $C.$
Thus
\begin{eqnarray*}
 \mathbb{P}\Big( \Big| X_r(u)-X_r(u_1)\Big|\geq \lambda,\ \Big|X_r(u_2)-X_r(u)\Big|\geq \lambda \Big)
&\leq&\frac{4C_1}{\lambda^{\gamma} } \Big[ u_2   - u_1  \Big]^\eta.
\end{eqnarray*}
Hence, by (\ref{ineq38}), if $\eta>1$, then $X(t)$  is $1/\alpha(x)-$strongly localisable
at $x$ with strong local form $X_x'(t)$.

If $X(t)$ has independent increments, then
\begin{eqnarray}
&&  \mathbb{P}\Big( \Big| X_r(u)-X_r(u_1)\Big|\geq \lambda,\ \Big|X_r(u_2)-X_r(u)\Big|\geq \lambda \Big)  \nonumber\\
&&\ \ \ \ \ \ \ \ \ \  \ =\ \mathbb{P}\Big( \Big| X_r(u)-X_r(u_1)\Big|\geq \lambda \Big) \, \mathbb{P}\Big( \Big|X_r(u_2)-X_r(u)\Big|\geq \lambda \Big).
\end{eqnarray}
By an argument similar to (\ref{sfnmbs}), it follows that
\begin{eqnarray*}
&& \mathbb{P}\Big( \Big| X_r(u)-X_r(u_1)\Big|\geq \lambda \Big)\  \leq \ \frac{4C_1 }{\lambda^{\gamma} } \Big[ u   - u_1  \Big]^\eta
\end{eqnarray*}
and
\begin{eqnarray*}
&& \mathbb{P}\Big( \Big|X_r(u_2)-X_r(u)\Big|\geq \lambda \Big) \  \leq \ \frac{4C_1 }{\lambda^{\gamma} } \Big[ u_2   - u   \Big]^\eta .
\end{eqnarray*}
Using the inequality $xy \leq (x+y)^2/4,$  $x, y \geq0,$   we have
\begin{eqnarray}
\mathbb{P}\Big( \Big| X_r(u)-X_r(u_1)\Big|\geq \lambda,\ \Big|X_r(u_2)-X_r(u)\Big|\geq \lambda \Big) \leq \frac{16C_1^2}{\lambda^{2 \gamma } } \Big[  u_2 - u_1 \Big]^{2\eta}.
\end{eqnarray}
Thus, if $2\eta> 1,$ by (\ref{ineq38}), then $X(t)$  is $1/\alpha(x)-$strongly localisable at $x$.
This completes the proof of the theorem.\qed

As an example to illustrate Theorem \ref{scrend}, consider the weighted MsLM.
Falconer  and  Liu  have proved that the weighted MsLM are localisable.
The following theorem shows that the weighted MsLM are not only localisable but also strongly localisable.
In particular, it shows that  the independent-increments  MsLM
are strongly localisable.

\begin{theorem}\label{themni}   Assume that the function $\alpha(u), u  \in [0,1],$ satisfies condition
 (\ref{coalph}). Let
 \[
Y(t)=\int_0^1 w(x)\mathbf{1}_{[0,\ t]}(x)M_\alpha ( dx), \ \ \ \ \  t \in [0, 1],
\]
 be a weighted multistable L\'{e}vy motion, where the function $ w(x), x \in [0, 1],$ is continuous.
Then   $Y(t)$ is $1/\alpha(x)-$strongly localisable
at all $x\in [0, 1]$ with strong local form $w(x)L_{\alpha (x) }(\cdot)$.
In particular, this implies that  $L_{I }(t)$
is $1/\alpha(x)-$strongly localisable
at all $x\in [0, 1]$ with strong local form $L_{\alpha (x) }(\cdot)$, an $\alpha (x)-$stable L\'{e}vy motion.
\end{theorem}
\noindent\emph{Proof.} It is known that $Y(t)$ is $1/\alpha(x)-$localisable
at all $x$ with strong local form $w(x)L_{\alpha (x) }(\cdot)$; see Falconer  and  Liu \cite{FL12}.
Set $f(t,x)=w(x)\mathbf{1}_{[0,\ t]}(x), \ t,x \in [0, 1].$
By (\ref{sdfdsdf}), the integrand of  $Y(t)$ satisfies, for all $t, v \in [0, 1]$ such that $v \leq t,$
\begin{eqnarray}
  \int_0^1 \Bigg| \frac{f(x+rt,s) - f(x+rv,s)  }{ r^{1/\alpha(x)} } \Bigg|^{\alpha(s)} ds
&=&\int_0^1 \Bigg| \frac{w(s)\mathbf{1}_{[x+rv, \, x+rt]}(s) }{ r^{1/\alpha(x)} } \Bigg|^{\alpha(s)} ds   \nonumber \\
&=&\int \Big| w(x+rz) \mathbf{1}_{[v,\  t ] }(z)\Big|^{\alpha (x+rz)} r^{(\alpha (x)-\alpha (x+rz))/\alpha (x)} dz \nonumber \\ &\leq&C_w\int   \mathbf{1}_{[v,\  t ] }(z)\,  r^{(\alpha (x)-\alpha (x+rz))/\alpha (x)} dz \nonumber \\
&\leq& 2C_w \ (t- v),  \nonumber
\end{eqnarray}
for all sufficiently small $r>0,$
where $s=x+rz$ and    $C_\omega = \sup_{z \in [0, 1]}|w(z)|^{\alpha(z)}. $
By the fact that the integrals of functions with disjoint supports are independent, it is easy to see that $Y(t)$ has independent increments, the first claim of the theorem follows by Theorem \ref{scrend}. In particular, since  $L_{I }(t) =\int_0^1  \mathbf{1}_{[0,\ t]}(x)M_\alpha (dx),$ the first claim of the theorem implies the second one with $w(x)=1, x \in [0, 1].$ \qed
\begin{remark}
By inspecting the proof of Falconer  and  Liu \cite{FL12}, we can see that $Y(t)$ is also $1/\alpha(x)-$localisable
at all $x$ with strong local form $w(x)L_{\alpha (x) }(t)$ when the function $ w(x), x \in [0, 1],$ is c\`{a}dl\`{a}g.
Hence, Theorem \ref{themni} holds true when the function $ w(x), x \in [0, 1],$ is c\`{a}dl\`{a}g.
% It is worth noting that for the real-world examples, the function $w(x)$ is c\`{a}dl\`{a}g and not continuous.
\end{remark}

\end{document}